\documentclass[12pt,reqno,a4paper]{amsart}
\usepackage{amsmath,amssymb,amsfonts,amsthm}
\usepackage[mathscr]{eucal}
\usepackage{hyperref}

\author{Hovhannes~M.~Khudaverdian}
\thanks{}
\address{School of Mathematics, University of Manchester, Sackville Street, Manchester, M60 1QD, United Kingdom}
\email{khudian@manchester.ac.uk, theodore.voronov@manchester.ac.uk}

\author{Theodore~Th.~Voronov}

\title[Differential forms and odd symplectic geometry]{Differential forms and odd symplectic geometry}

\newtheorem{thm}{Theorem}

\newtheorem{lm}{Lemma}[section]

\newtheorem{cor}{Corollary}[section]

\theoremstyle{definition}
\newtheorem{rem}{Remark}[section]

\def\co{\colon\thinspace}

\renewcommand{\i}{\mathrm{i}}

\newcommand{\schl}{{[\![}}
\newcommand{\schr}{{]\!]}}

\DeclareMathOperator{\Can}{Can}

\DeclareMathOperator{\Ber}{Ber} \DeclareMathOperator{\Ker}{Ker}
\DeclareMathOperator{\Def}{Def} \DeclareMathOperator{\Ind}{Ind}
\renewcommand{\Im}{\mathop{\mathrm{Im}}}

 \renewcommand{\P}{\Pi}

 \DeclareMathOperator{\const}{const}

\DeclareMathOperator{\Ann}{Ann} \DeclareMathOperator{\Sp}{Sp}

\renewcommand{\div}{\mathop{\mathrm{div}}}

\newcommand{\der}[2]{{\frac{\partial {#1}}{\partial {#2}}}}

\newcommand{\dder}[3]{{\frac{\partial^2 {#1}}{\partial {#2}\partial {#3}}}}

\newcommand{\Z}{{\mathbb Z_{2}}}
\newcommand{\ZZ}{{\mathbb Z}}
\newcommand{\p}{\partial}

\newcommand{\widebar}{\overline}

\renewcommand{\a}{\alpha}
\renewcommand{\b}{\beta}

\def\s{\sigma}

\def\O{\Omega}
\def\D{\Delta}

\def\o{\omega}

\renewcommand{\r}{{\rho}}

\newcommand{\x}{{\xi}}
\def\d{\delta}

\begin{document}
\begin{abstract}
We recall the main facts about the odd Laplacian acting on
half-densities on an odd symplectic manifold and discuss a
homological interpretation for it suggested recently by P.~{\v
S}evera. We study the relationship of odd symplectic geometry with
classical objects. We show that the Berezinian of a canonical
transformation for an odd symplectic form is a polynomial in matrix
entries and a complete square. This is a simple but fundamental
fact, parallel to Liouville's theorem for an even symplectic
structure. We draw attention to the fact that the de Rham complex on
$M$ naturally admits an action of the supergroup of all canonical
transformations of $\Pi T^*M$. The infinitesimal generators of this
action turn out to be the classical `Lie derivatives of differential
forms along multivector fields'.
\end{abstract}
\maketitle

Odd symplectic geometry (more generally, odd Poisson geometry) or
the geometry of odd brackets is the mathematical basis of the
Batalin--Vilkovisky method~\cite{bv:perv,bv:vtor,bv:closure} in
quantum field theory.

Odd symplectic geometry possesses features connecting it with both
classical (``even'') symplectic geometry and Riemannian geometry. In
particular, \textit{odd Laplace operators} arise naturally on an odd
symplectic manifold, i.e., the second order differential operators
whose principal symbol is the odd quadratic form corresponding to
the odd bracket~\cite{hov:deltabest}.  The key difference from the
Riemannian case is that the definition of an odd Laplace operator,
in general, requires an extra piece of data besides the ``metric'',
namely, a choice of a volume form (even for a Laplacian acting on
functions). This is due to the fundamental fact that on an odd
symplectic manifold there is no invariant volume
element~\cite{hov:deltabest}.

However, as   discovered by one of the authors, there is one
isolated case where an odd Laplacian is defined canonically by the
symplectic structure without any extra
data~\cite{hov:max,hov:semi,hov:proclms}. It is an operator acting
on densities of weight $1/2$ (half-densities or semidensities). This
fact is not obvious, and there is no simple explanation. A known
proof is based on an analysis of the canonical transformations of
the odd bracket. In works~\cite{tv:laplace1, tv:laplace2,
tv:laplace2bis} further phenomena related with odd Laplacians on odd
Poisson manifolds were discovered, such as the existence of a
natural `master' groupoid acting on volume forms,   its orbits
corresponding to Laplacians on half-densities.  The symplectic case
is distinguished by the existence of a distinguished orbit, which
gives the ``canonical'' operator.

In a very interesting recent paper~\cite{severa:originbv}, P.~{\v
S}evera suggested a homological interpretation of the canonical odd
Laplace operator on half-densities  as one of the higher
differentials in a certain natural spectral sequence associated with
the odd symplectic structure.

In our paper we  in particular  discuss this interpretation and show
that there is a   simple but fundamental underlying fact from linear
algebra, concerning the Berezinian of a canonical transformation for
an odd symplectic bracket. It is the formula
\begin{equation}\label{eq:berj1}
    \Ber J =(\det J_{00})^2
\end{equation}
for $J$ in the  odd  symplectic supergroup, where $J_{00}$ is the
even-even block. Hence the Berezinian is an entire rational function
and, moreover, a complete square. There are many geometric facts
related with formula~\eqref{eq:berj1}, which  can be found in the
literature on odd  brackets and the BV formalism. See, for example,
~\cite{ass:bv, ass:symmetry,  hov:deltabest, hov:bv, hov:semi}. We
want to draw attention to it as   a simple identity for matrices. In
view of it, half-densities on an odd symplectic manifold are
`tensor' objects, i.e., transforming according to a polynomial
representation. They can be seen as \textit{virtual differential
forms} on a Lagrangian surface. When such a surface is fixed, they
become (isomorphic to) actual forms. More precisely, for a manifold
or supermanifold $M$, we can identify (pseudo)integral forms on $M$,
i.e., multivector densities, with half-densities on the odd
symplectic manifold $\Pi T^*M$. (Pseudo)integral forms are related
with (pseudo)differential forms by a sort of `Fourier transform'.
Therefore we see that in the space of differential forms on an
ordinary manifold, there is a natural representation of the
supergroup of canonical transformations of the odd bracket. We give
a clear description of this action in classical terms.  The
invariance of the de Rham differential under such a supergroup,
which is absolutely transparent, is equivalent to the existence of
the canonical odd Laplacian, but expressed in a different language.

\section{Recollection of the canonical odd Laplacian}

In this section we review the construction of the odd Laplacian on
half-densities due to~\cite{hov:max}. See also~\cite{hov:semi,
hov:proclms} and \cite{tv:laplace1}.

Let $M$ be a supermanifold endowed with an odd symplectic structure,
given by an odd $2$-form $\o$. We shall refer to such supermanifolds
as  \textit{odd symplectic manifolds}. (We always skip the prefix
`super-' unless required to avoid confusion.) Later we shall discuss
the more general case of an odd Poisson manifold. A brief definition
of the odd Laplacian acting on half-densities on $M$  follows.

Consider a cover of $M$ by Darboux charts, in which the symplectic
form takes the canonical expression $\o=dx^id\x_i$. Here $x^i$,
$\x_i$ are canonically conjugate variables of opposite parity. We
assume that the $x^i$ are even; hence the $\x_i$, odd. Let $Dy$, for
any kind of variables $y$, stand for the Berezin volume element.
Then half-densities on $M$ locally look like
$\s=s(x,\x)(D(x,\x))^{1/2}$. (Notice that we skip  questions related
with orientation.) We set
\begin{equation}\label{eq:delta}
    \D \s:=\dder{s}{x^i}{\x_i}\,\bigl(D(x,\x)\bigr)^{1/2},
\end{equation}
in Darboux coordinates, and call $\D$, the \textit{canonical odd
Laplacian on half-densities}.

The simplicity of formula~\eqref{eq:delta} is very deceptive. The
expression $\der{}{x^i}\der{f}{\x_i}$ was originally suggested by
Batalin and Vilkovisky, and is the famous `BV operator'. However,
the trouble is that it is not well-defined on functions (actually,
on any objects) unless we fix a volume form, which should therefore
enter the definition.  The geometrically invariant construction for
functions, using a volume form, was first given
in~\cite{hov:deltabest}. \textit{There is no canonical volume form
on an odd symplectic manifold} (unlike even symplectic manifolds,
enjoying the Liouville form). In particular, the coordinate volume
form $D(x,\x)$ for Darboux coordinates is \textbf{not} preserved by
the (canonical) coordinate transformations (see later). Hence the
invariance of the operator $\D$ given by~\eqref{eq:delta} is a deep
geometric fact.

As we showed in~\cite{tv:laplace1}, on any odd Poisson, in
particular, odd symplectic,  manifold there is a natural
\textit{master groupoid} of `changes of volume forms' $\r\mapsto
e^S\r$ satisfying the master equation $\D_{\r}e^{S/2}=0$ (note $1/2$
in the exponent; without it there would be no groupoid). Here
$\D_{\r}$ is the odd Laplacian on functions with respect to the
given volume form $\r$. It is defined by $\D_{\r}f:=\div_{\r}X_f$,
where $X_f$ is the Hamiltonian vector field corresponding to $f$.
(See~\cite{hov:deltabest}; note also~\cite{yvette:divergence} for
another approach.) In a similar way one can define the odd Laplacian
on any densities
--- again, depending on a chosen volume form.  Now, half-densities are
distinguished from densities of other weights precisely by the fact
that for them the corresponding odd Laplacian would depend only on
the orbit of a volume form with respect to the action of the above
groupoid~\cite{tv:laplace1}. It turns out that on an odd symplectic
manifold, all Darboux coordinate volume forms belong to the same
orbit of the master groupoid. We can regard it as a `preferred
orbit'; hence, in the absence of an invariant volume form, the odd
Laplacian on half-densities defined by an arbitrary  Darboux
coordinate volume form  is invariant. It is just~\eqref{eq:delta}.

\section{Homological interpretation of the odd Laplacian}

Now we are going to approach $\D$ on half-densities from a very
different angle.

Let $\O(M)$ be the space of all pseudodifferential forms on $M$,
i.e., functions on $\P TM$. (As usual, $\P$ stands for the parity
reversion functor on vector spaces, vector bundles, etc.) In
coordinates such functions have the form  $s=s(x,\x,dx,d\x)$, where
the differentials of coordinates are commuting variables of parity
opposite to that of the respective coordinate. In our case $dx^i$
are odd and $d\x_i$ are even. We do not assume that functions
$s(x,\x,dx,d\x)$ are polynomial in $d\x_i$.  Of course they are
(Grassmann) polynomial in $dx^i$, because these variables are odd.

Consider the odd symplectic form $\o$. Since $\o^2=0$,
multiplication by $\o$ can be considered as a differential. Define
the operator $D=d+\o$, where $d$ is the de Rham differential. Since
$d\o=0$, it follows that $D^2=0$ and we have a `double complex'
$\bigl(\O(M),D=\o+d\bigr)$. \textbf{Warning:} here a
\textit{complex} means just a $\Z$-graded object.

The reader should bear in mind that since $\o=d\Theta$ for some even
$1$-form $\Theta$, which is true globally, we have
$D=e^{-\Theta}\circ d \circ e^{\Theta}$ and the multiplication by
the inhomogeneous differential form $e^{\Theta}$ sets an isomorphism
between the complexes $\bigl(\O(M), D\bigr)$ and $\bigl(\O(M),
d\bigr)$. It follows that $H(\O(M), D)$ is isomorphic to $H(\O(M),
d)$, which is just the de Rham cohomology of the underlying manifold
$M_0$. (Note that the isomorphism $e^{\Theta}$ preserves only
parity, but not $\ZZ$-grading, even if we restrict it to
{differential} forms on $M$, i.e., polynomials in $dx,d\x$.)

The operator $D=\o+d$ was introduced in~\cite{severa:originbv}. The
idea was to consider the spectral sequence for $\bigl(\O(M),D\bigr)$
regarded as a double complex. We shall follow it in a form best
suiting our purposes and which is slightly different
from~\cite{severa:originbv}. (In particular, we do not assume
grading in the space of forms.)

Although there is no $\ZZ$-grading present, single or double, one
can still develop the machinery of spectral sequences as follows.

We define linear relations (see~\cite{maclane:homology}) on $\O(M)$:
\begin{equation*}
    \p_0:=\bigl\{(\a,\b) \in \O(M)\times \O(M)\ |\  \o\a=\b\bigr\},
\end{equation*}
and
\begin{multline*}
    \p_r:=\bigl\{(\a,\b) \in \O(M)\times \O(M)\ |\  \exists \a_1,\ldots, \a_r\in \O(M) \ :
    \ \o\a=0, \\
    d\a+ \o\a_1=0, \ \ldots, 
    \ d\a_{r-2}+ \o\a_{r-1}=0, \ d\a_{r-1}+ \o\a_{r}=\b\bigr\}
\end{multline*}
for all $r=1,2, 3, \ldots \ {}$  We also set $\p_{-1}:=\{(\a,0)\}$.
We have subspaces $\Ker\p_r$, $\Def\p_r$ (the \textit{domain of
definition}), $\Ind\p_r$ (the \textit{indeterminancy}), and
$\Im\p_r$ in $\O(M)$, and by a direct check
\begin{align*}
    \Im \p_r&\subset \Ker \p_r\,,\\
\Def \p_r&=\Ker\p_{r-1}\,,\\
\Ind \p_r&=\Im\p_{r-1}\,.
\end{align*}
That is, we have a sequence of differential relations on $\O(M)$,
defining a spectral sequence $(E_r,d_r)$  where
\begin{equation*}
    E_{r}:=\frac{\Ker \p_{r-1}}{\Im\p_{r-1}}=\frac{\Def \p_{r}}{\Ind\p_{r}}
\end{equation*}
and the homomorphism $d_r\co E_r\to E_r$ is induced by $\p_r$ in the
obvious way. (In fact, differential relations like this is the
shortest way of defining spectral sequences,
see~\cite[p.~340]{maclane:homology}.)

Clearly $E_0=\O(M)$. The relation $\p_0$ is simply the graph of the
linear map $d_0\co \O(M)\to \O(M)$, $d_0\a=\o\a$. What is $E_1$?

\begin{thm} The space $E_1$ can be naturally identified with the space
of half-densities on $M$.
\end{thm}

A proof consists of two independent steps. First, we find the
cohomology of $d_0$ using algebra. Second, we identify the result
with a geometrical object. The first part goes as follows.

The operator $d_0=\o$ is a Koszul type differential, since in an
arbitrary Darboux chart $\o=dx^id\x_i$. Introduce a $\ZZ$-grading by
the degree in the odd variables $dx^i$. The operator $d_0$ increases
the degree by one. (This grading is \textbf{not} preserved by
changes of coordinates.) From general theory it follows that the
cohomology should be concentrated in the ``maximal degree''. Indeed,
suppose that $\dim M=n|n$ and consider the linear operator $H$ on
pseudodifferential forms defined as follows. For
$\s=\s(x,\x,dx,d\x)$,
\begin{equation*}
    H\s(x,\x,dx,d\x):=\int_0^1 \!dt\,t^{n-1}\, \dder{\s}{dx^i}{d{\x}_i}
    (x,\x,t^{-1}dx,t\,d\x)\,,
\end{equation*}
--- notice the similarity with the $\D$-operator. The operator $H$ is
well defined on all forms of degree  less than $n$ in $dx^i$ and on
forms of `top' degree if they vanish at $d\x_i=0$. (In both cases
there will be no problem with division by $t$.) For forms on which
$H$ makes sense one can check that
\begin{equation*}
    (Hd_0+d_0H)\s=\s\,.
\end{equation*}
In particular, if a form $\s$ is $d_0$-closed and of degree less
than $n$ in $dx^i$, then $\s=d_0 H\s$. The same applies for a top
degree form taking a non-zero value at $d\x_i=0$. Hence the
$d_0$-cohomology ``sits on'' pseudodifferential forms of degree $n$
in $dx^i$ that do not depend on $d\x_i$:
\begin{equation*}
    \s=s(x,\x)\,dx^1\ldots dx^n.
\end{equation*}
No non-zero form of this appearance can be cohomologous to zero:
indeed, any $d_0$-exact form, $d_0\tau=\o\tau$, vanishes at
$d\x_i=0$.

\textit{Hence, each $d_0$-cohomology class has a unique
representative in a given Darboux coordinate system $x^i,\x_i$. It
is obtained by taking an arbitrary form from the class, extracting
its component of degree $n$ in $dx^i$ and evaluating  at $d\x_i=0$.}
By applying this to the class of $dx^1\ldots dx^n$, we immediately
arrive at
\begin{lm} \label{lem1} Elements of the cohomology space $E_1=H(\O(M),\o)$
are represented in Darboux coordinates as classes
\begin{equation*}
    \s=s(x,\x)\,[dx^1\ldots dx^n],
\end{equation*}
where under a change of Darboux coordinates
\begin{align*}
    x^i&=x^i(x',\x'),\\
\x_i&=\x_i(x',\x')
\end{align*}
the class $[dx^1\ldots dx^n]$ transforms as follows:
\begin{equation*}
    [dx^1\ldots dx^n]=\det J_{00}\cdot[dx^{1'}\ldots
    dx^{n'}]\,.
\end{equation*}
Here $J_{00}=\der{x}{x'}$ is the even-even block of the Jacobi
matrix $J=\der{(x,\x)}{(x',\x')}$. 
\end{lm}

To  better appreciate the statement, notice that
\begin{equation*}
    dx^i=dx^{i'}\der{x^i}{x^{i'}}+d\x_{i'}\der{x^i}{\x_{i'}}.
\end{equation*}
Hence
\begin{equation*}
    dx^1\ldots dx^n=dx^{1'}\ldots
    dx^{n'}\cdot \det\left(\der{x^i}{x^{i'}}\right)+ \text{terms containing
    $d\x_{i'}$}\,.
\end{equation*}
Passing to cohomology is equivalent to discarding these lower order
terms.

What kind of geometrical object is this?

\begin{lm} \label{lem2} Objects of the form
    $\s=s(x,\x)\,[dx^1\ldots dx^n]$,
in Darboux coordinates, with the transformation law given in
Lemma~\ref{lem1} can be identified with half-densities on $M$.
\end{lm}

This is the crucial claim. There is a simple but fundamental fact
from linear algebra behind Lemma~\ref{lem2}, which will be proved in
the next section.

The transformation law for $[dx^1\ldots dx^n]$ can be obtained from
the formal ``law'' $[dx^i]=[dx^{i'}]\der{x^i}{x^{i'}}$.
Unfortunately, it does not define a geometric object, because it
does not obey the cocycle condition. In a way, it is only a
`virtual' transformation law, which will make sense only if an extra
structure is imposed on $M$.

Now as we have the space $E_1$, let us check the differential $d_1$
on it. It is induced  by the differential relation $\p_1$ on
$\O(M)$. Take an element $\s=s(x,\x)\,[dx^1\ldots dx^n]\in E_1$,
take its representative $\a=s(x,\x)\, dx^1\ldots dx^n$ and consider
$\b\in\O(M)$  such that $d\a+\o\a_1=\b$, for $\a_1\in \O(M)$. We
will have $[\b]=d_1\s$ for the class $[\b]$ in $E_1$. Notice that
$d\a=d\x_i\der{s}{\x_i}\, dx^1\ldots dx^n$ and it will vanish at
$d\x_i=0$, therefore it is an $\o$-exact form, according to our
previous analysis. \textit{Thus $d_1=0$ identically and $E_2=E_1$.}

Consider $d_2$ on $E_1=E_2=H(\O(M),\o)$. By definition, $d_2$ maps
the class $\s=s(x,\x)\,[dx^1\ldots dx^n]$, with a local
representative $\a=s(x,\x)\, dx^1\ldots dx^n$, to the class of
$\b\in \O(M)$ such that $d\a+\o\a_1=0$, $d\a_1+\o\a_2=\b$, for some
$\a_1$ and $\a_2$.  We may set $\a_1:=-Hd\a$, where $H$ is the
homotopy operator defined above, and  $\b:=d\a_1=-dHd\a$. Directly:
\begin{multline*}
    Hd\a=H\Bigl(d\x_i\der{s}{\x_i}\, dx^1\ldots dx^n\Bigr)=\sum (-1)^{i +\tilde
s} \der{s}{\x_i}\, dx^1\ldots \widehat{dx^i} \ldots dx^n
\end{multline*}
and
\begin{multline*}
    \b=-dHd\a=-d \sum (-1)^{i-1+\tilde
s} \der{s}{\x_i}\, dx^1\ldots \widehat{dx^i} \ldots dx^n=\\
-dx^j\der{}{x^j}\sum (-1)^{i +\tilde s} \der{s}{\x_i}\, dx^1\ldots
\widehat{dx^i} \ldots dx^n + \text{\ lower order terms in $dx$}=\\
-\dder{s}{x^i}{\x_i}\,dx^1\ldots dx^n + \text{\ lower order terms in
$dx$}\,.
\end{multline*}
Hence in $E_1$ we get:
\begin{equation*}
    d_2\s=d_2 \bigl(s(x,\x)\,[dx^1\ldots dx^n]\bigr)=-\dder{s}{x^i}{\x_i}\,[dx^1\ldots
    dx^n]=-\D\s\,,
\end{equation*}
which is quite remarkable. What about the space $E_3$ and the
differential $d_3$, and so on?

It is not hard to notice that the cohomology of the $\D$-operator on
half-densities on $M$ is isomorphic to the de Rham cohomology of the
underlying ordinary manifold $M_0$ (we shall say more about this
later). Locally the cohomology vanishes except for constants:
$\s=\const\cdot [dx^1\ldots dx^n]$. Thus $d_3=0$, and $E_4=E_3$; the
same continues for $d_4=0$,  $E_5=E_4=E_3$, and so on. We arrive at
the following statement (which was the main result
of~\cite{severa:originbv}):

\begin{thm} \label{thm2} With the identification of the space $E_1=H(\O(M),\o)$
with half-densities on $M$, the differential $d_1$ vanishes and the
next differential $d_2$ coincides up to a sign with the canonical
odd Laplacian. The spectral sequence $(E_r,d_r)$ degenerates at the
term $E_3$, which is the cohomology of the operator $\D$.
\end{thm}

The importance of Theorem~\ref{thm2}  is in the fact that it gives
an alternative proof of the invariance of the odd Laplacian on
half-densities $\D$, by identifying  it with an operator in a
spectral sequence invariantly associated with the odd symplectic
structure.

\section{Berezinian of a canonical transformation}

Consider a vector space $V=V_0\oplus V_1$ with an odd symplectic
structure, i.e., an odd non-degenerate antisymmetric bilinear form.
(A choice of `antisymmetric' or `symmetric' does not make any
difference.) Necessarily $\dim V=n|n$. We call matrices preserving
this form, \textit{symplectic}. This should not cause problems; when
comparing them with ordinary symplectic matrices  corresponding to
an even symplectic structure, we shall make the reference to the
parity of the bilinear form explicit.

\begin{thm} \label{thm:block} Suppose that $J$ is a symplectic matrix for an odd symplectic space.
Let
\begin{equation*}
    J=\begin{pmatrix}
        J_{00} & J_{01} \\
        J_{10} & J_{11} \\
      \end{pmatrix}
\end{equation*}
be its standard block decomposition. Then
\begin{equation}
    \Ber J =(\det J_{00})^2\,.
\end{equation}
\end{thm}

\begin{proof} We can write the matrix of our symplectic form as
\begin{equation*}
    B=\begin{pmatrix}
        0 & 1 \\
        1 & 0 \\
      \end{pmatrix}\,.
\end{equation*}
The relation for $J$ is $JBJ\,^T=B$, where the operation of matrix
transpose takes into account the parities of the blocks:
\begin{multline*}
    \begin{pmatrix}
        J_{00} & J_{01} \\
        J_{10} & J_{11} \\
      \end{pmatrix}
      \begin{pmatrix}
        0 & 1 \\
        1 & 0 \\
      \end{pmatrix}
      \begin{pmatrix}
        J_{00} & J_{01} \\
        J_{10} & J_{11} \\
      \end{pmatrix}^T=
\begin{pmatrix}
        J_{00} & J_{01} \\
        J_{10} & J_{11} \\
      \end{pmatrix}
      \begin{pmatrix}
        0 & 1 \\
        1 & 0 \\
      \end{pmatrix}
      \begin{pmatrix}
        J_{00}^{\,T} & J_{10}^{\,T} \\
        -J_{01}^{\,T} & J_{11}^{\,T} \\
      \end{pmatrix}.
\end{multline*}
Hence we obtain
\begin{align}
    J_{00}J_{01}^{\,T}=\left(J_{00}J_{01}^{\,T}\right)^T, \label{eq:symp1}\\
    J_{11}J_{10}^{\,T}=-\left(J_{11}J_{10}^{\,T}\right)^T, \label{eq:symp2}\\
    J_{00}J_{11}^{\,T}+J_{01}J_{10}^{\,T}=1 \label{eq:symp3}.
\end{align}
From~\eqref{eq:symp3} we may express
\begin{multline*}
    J_{11}={J_{00}^{\,T}}^{-1}+J_{10}{J_{01}^{\,T}}{J_{00}^{\,T}}^{-1} \\
    ={J_{00}^{\,T}}^{-1}+J_{10}J_{00}^{-1}J_{01}, \text{\quad   taking into
    account~\eqref{eq:symp1}}.
\end{multline*}
We arrive at the identity
\begin{equation}
    J_{11}-J_{10}J_{00}^{-1}J_{01}={J_{00}^{\,T}}^{-1}.\label{eq:id}
\end{equation}
Therefore
\begin{equation*}
    \Ber J=\frac{\det
    J_{00}}{\det\bigl(J_{11}-J_{10}J_{00}^{-1}J_{01}\bigr)}=
    \frac{\det
    J_{00}}{\det{J_{00}^{\,T}}^{-1}}=\left(\det
    J_{00}\right)^2.
\end{equation*}
\end{proof}

Notice that the steps we followed in the proof are the same as may
be used for proving the classical Liouville theorem, i.e., that
$\det J=1$ for an ordinary symplectic matrix $J\in \Sp(2n)$. The
decomposition of $V$ in that case will be the decomposition into the
sum of Lagrangian subspaces and an identity similar to~\eqref{eq:id}
will be valid. The difference will arise only when calculating the
determinant: instead of the ratio  of the determinants of the
blocks, there will the product, which will give $1$ instead of $\det
J_{00}^2$.

It is easy to generalize. Let $V$ be a vector (super)space with a
symplectic structure, even or odd. Consider its decomposition into
the sum of two Lagrangian subspaces. In the even case they will have
the same dimensions; in the odd case, the opposite, i.e., $p|{n-p}$
and ${n-p}|p$. Denote the chosen decomposition by $V=V_0\oplus V_1$.
Here the indices have nothing to do with parity. By picking
`canonically conjugate' bases in $V_0$ and $V_1$ we arrive at a
picture formally the same as above. The Berezinian can be calculated
using the corresponding block decompositions. It will be either
\begin{equation*}
    \Ber J=\Ber J_{00}\cdot \Ber \bigl(J_{11}-J_{10}J_{00}^{-1}J_{01}\bigr)
\end{equation*}
or
\begin{equation*}
    \Ber J=\frac{\Ber J_{00}}{\Ber
    \bigl(J_{11}-J_{10}J_{00}^{-1}J_{01}\bigr)}
\end{equation*}
depending on the parity of the symplectic form (in the odd case the
`formats' of the matrix blocks will be the opposite, hence
division). Then the analog of the identity~\eqref{eq:id} should be
applied. For even symplectic structure we thus obtain the analog of
Liouville's theorem, and for odd, we arrive at

\begin{thm} \label{thm:lagrange} Let $J$ be a symplectic transformation of an odd symplectic space $V$.
Then for an arbitrary decomposition into the sum of Lagrangian
subspaces, $V=V_0\oplus V_1$ \emph{(indices not indicating parity)},
the identity
\begin{equation}\label{eq:berj}
    \Ber J=\left(\Ber J_{00}\right)^2
\end{equation}
holds.
\end{thm}

\begin{rem} Theorem~\ref{thm:block}  gives, in particular, that the
Berezinian of a symplectic matrix is a polynomial in the matrix
entries and, moreover, a complete square. This is somewhat masked in
the more general Theorem~\ref{thm:lagrange}.
\end{rem}

There is an `abstract' argument parallel to the calculation above.
Consider a decomposition into Lagrangian subspaces $V=V_0\oplus
V_1$, for an even or odd symplectic form. (This works in the same
way for symmetric forms.) Consider the dual space $V_1^*$. We have
$V_1^*\cong \Ann V_0\subset V^*$ and, using the form, can replace
the annihilator by the orthogonal complement:
\begin{equation*}
    V_1^*\cong V_0^{\perp}\subset V
\end{equation*}
for the even form or
\begin{equation*}
    V_1^*\cong \Pi V_0^{\perp}\subset \Pi V
\end{equation*}
for the odd form. Recalling that $V_0$ and $V_1$ are Lagrangian
subspaces, we get $V_1^*\cong V_0$ or $V_1^*\cong \Pi V_0$. Hence
\begin{align*}
    \Ber V=\Ber (V_0\oplus V_1)&=
    \Ber V_0\otimes \Ber V_1\\
    &=\Ber
    V_0\otimes \Ber V_0^*=1
    \intertext{(even symplectic form) or}
&=\Ber V_0\otimes \Ber \Pi V_0^*=(\Ber V_0)^{\otimes 2}
\end{align*}
(odd symplectic form). Equalities here mean  natural isomorphisms.

A different abstract argument based on the well-known interpretation
of the space $\Ber V$ as the cohomology of a Koszul complex and
justifying the equality $\Ber V=(\Ber V_0)^{\otimes 2}$ for an odd
symplectic space was given in~\cite{severa:originbv}.

A weak point of abstract arguments is that they do not really give
information about matrices, which is necessary in applications such
as a proof of Lemma~\ref{lem2}.

\section{``Supersymmetries'' of differential forms}

In this section we change   viewpoint. We would like to phrase the
previous constructions entirely in the language of `classical'
differential-geometric objects. In this way we shall see how the
canonical odd Laplacian on half-densities on an odd symplectic
manifold can be seen as a `classical' object equipped with extra
symmetries. 

Let $M$ now stand for an arbitrary manifold or supermanifold.
Previously we worked with   odd symplectic manifolds. It is known
that any such odd symplectic manifold can be  non-canonically
identified with $\Pi T^*M$ considered with the natural odd bracket,
for some $M$. A change of   identifying symplectomorphism is
equivalent to a symplectomorphism (or canonical transformation) of
the space $\Pi T^*M$. Therefore, we can restrict ourselves to
objects on $\Pi T^*M$, but should analyze them from   the viewpoint
of the larger supergroup of all canonical transformations of $\Pi
T^*M$, not just that of diffeomorphisms of $M$.

We shall consider multivector fields (and multivector densities)
and differential forms on $M$. When $M$ is a supermanifold, we
actually speak of pseudodifferential forms.

Multivector fields on $M$ are identified with functions on $\Pi
T^*M$. In local coordinates, we have $X=X(x,x^*)$, where $x^a$ are
coordinates on $M$ and $x^*_a$ are the corresponding coordinates on
the fibers, of the opposite parity, transforming as
$x^*_a=\der{x^{a'}}{x^a}\,x^*_{a'}$. There is no problem with
canonical transformations of $\Pi T^*M$ acting on multivector fields
on $M$ --- it is just the pull-back of functions.

Multivector densities on $M$ have the form $\s=s(x,x^*)\,Dx$ and at
first glance it is not obvious how a transformation mixing $x$ and
$x^*$ can be applied to them. However, in view of
Theorems~\ref{thm:block} and \ref{thm:lagrange}, for  a canonical
transformation $F\co \Pi T^*M\to \Pi T^*M$ one can set
\begin{equation*}
    F^*\s:=s\bigl(y(x,x^*),y^*(x,x^*)\bigr)\,\Ber\left(\der{y}{x}\right)\,Dx
\end{equation*}
if  in coordinates $F\co (x,x^*)\mapsto \bigl(y=y(x,x^*),
y^*=y^*(x,x^*)\bigr)$. This is a well-defined action. In other
words, we identify multivector densities on $M$ with half-densities
on $\Pi T^*M$ and apply the natural  action, taking into account
identity~\eqref{eq:berj}. In integration theory, multivector
densities are known as integral forms (more precisely,
pseudointegral, if we insist on   differentiating between arbitrary
smooth functions and polynomials). Therefore we can make a  remark:
integral forms on an arbitrary supermanifold $M$ are the same  as
half-densities on the odd symplectic manifold $\Pi T^*M$. In this
language we see that integral forms have more symmetries than those
obvious ones  given by diffeomorphisms of $M$.

Consider  now pseudodifferential forms on $M$, i.e., functions on
$\Pi TM$. They are related with (pseudo)integral forms, i.e.,
multivector densities on $\Pi T^*M$ by the Fourier transform:
\begin{equation}\label{eq:fourier}
    \o(x,dx)=\int\limits_{\Pi T^*_xM}\! Dx^*\, e^{\i dx^ax^*_a}\,s(x,x^*)
\end{equation}
and conversely
\begin{equation}\label{eq:invfourier}
    s(x,x^*)=\const\int\limits_{\Pi T_xM}\! D(dx)\,
    e^{-\i dx^ax^*_a}\,\o(x,dx)\,.
\end{equation}
From here we obtain the action of the canonical transformations of
the odd symplectic manifold $\Pi T^*M$ on forms on $M$ as follows:
\begin{multline*}
    (F^*\o)(x,dx)=\const\iint_{\Pi T^*_xM\times \Pi T_yM} Dx^* D(dy)
    \,e^{\i \left(dx^ax^*_a-dy^ay^*_a(x,x^*)\right)}\cdot \\
    \o\bigl(y(x,x^*),dy\bigr)\,\Ber\der{y}{x}(x,x^*)\,,
\end{multline*}
where, as above, $F\co (x,x^*)\mapsto \bigl(y=y(x,x^*),
y^*=y^*(x,x^*)\bigr)$. In general, this action is non-local.

We shall consider the representation of the infinitesimal canonical
transformations of $\Pi T^*M$ on forms and multivector densities on
$M$. As it  turns out, the description in both cases will be very
simple.

For the odd symplectic manifold $\Pi T^*M$, the canonical odd
Laplacian on half-densities $\D$ on $\Pi T^*M$ is just the familiar
divergence of multivector densities $\d$ on $M$. Indeed,
\begin{equation*}
    \d\s= \dder{s}{x^a}{x^*_a} \,(x,x^*)\,Dx\,,
\end{equation*}
if $\s=s(x,x^*)\,Dx$. Consider the infinitesimal canonical
transformation of $\Pi T^*M$ generated by a function
(``Hamiltonian'') $H=H(x,x^*)$. Denote  by $L_H$ the corresponding
Lie derivative. Notice that from the viewpoint of $M$, the function
$H$ is a multivector field.
\begin{thm} On multivector densities (= integral forms) on $M$, the
Lie derivative w.r.t. the infinitesimal canonical transformation of
$\Pi T^*M$ generated by  $H$ is given by the formula:
\begin{equation*}
    L_H=[\d,H],
\end{equation*}
where at the r.h.s. stands the commutator of the divergence operator
$\d$ and  multiplication by the multivector field $H$.
\end{thm}

A proof can be given by a direct computation. It fact, the statement
mimics a similar and more general statement concerning odd Laplace
operators acting on densities of various weights,
see~\cite{tv:laplace1}.

\begin{cor} \label{cor:invar} The operator $\d$ on multivector densities on $M$ is
invariant under all canonical transformations of the odd symplectic
manifold $\Pi T^*M$. \emph{(At least those given by a Hamiltonian.)}
\end{cor}
\begin{proof} We need to show that $\d$ commutes with all Lie derivatives
$L_H$. Indeed, $[\d,L_H]=[\d,[\d,H]]=[\d^2,H]=0$, since $\d^2=0$.
\end{proof}

We can adopt the following viewpoint. Suppose we do not know
anything about the operator $\D$ on half-densities in odd symplectic
geometry. Instead we concentrate on a familiar object, the operator
$\d$ on multivector densities on a manifold $M$. The operator $\d$,
as shown, is invariant under much larger group of transformations
than just diffeomorphisms of $M$. It is invariant under
symplectomorphisms of $\Pi T^*M$. We can then take $\d$ as the
definition of $\D$ for $\Pi T^*M$. Since any odd symplectic manifold
$N$ is symplectomorphic to some $\Pi T^*M$, we can use this to
\textbf{define} $\D$ on $N$. The invariance of $\D=\d$ for $\Pi
T^*M$ under symplectomorphisms of $\Pi T^*M$ shows that $\D$ on $N$
is well-defined, i.e., its action on half-densities on $N$ does not
depend on an arbitrary choice of the identifying symplectomorphism
$N\cong \Pi T^*M$.

It may seem that there is a gap in such an argument as the
invariance was proved only infinitesimally or, equivalently,  for
transformations that can be included into a Hamiltonian flow. In
fact, there is no gap. Consider the supergroup $\Can \Pi T^*M$ of
all canonical transformations. If $M$ is an ordinary manifold, the
structure of this supergroup was described in~\cite{hov:semi}. It is
the product of the three subgroups:

\begin{enumerate}
  \item Transformations induced by diffeomorphisms of $M$;
  \item Shifts in the fibers of $\Pi T^*M$  of the form
\begin{equation*}
    F^*x^a=x^a, \quad
    F^*x^*_a=x^*_a+\der{\Phi}{x^a}\,,
\end{equation*}
where $\Phi=\Phi(x)$ is an odd function on $M$;
  \item Transformations identical on the submanifold $M\subset \Pi
  T^*M$.
\end{enumerate}

Since $\d$ is invariant under diffeomorphisms of $M$, all that
remains is to study transformations of types 2 and 3. Canonical
transformations of types 2 and 3  can be included into Hamiltonian
flows. Indeed, for type 2  one can take the flow with the
Hamiltonian $\Phi$. For type 3 there is also a Hamiltonian flow,
with the Hamiltonian of the form $\Psi=\Psi^{ab}(x,x^*)x^*_a x^*_b$,
as shown in~\cite{hov:semi}. Therefore,  transformations of types 2
and 3  are covered by the argument in the proof of
Corollary~\ref{cor:invar}, and this completes the proof.

Now let us turn our attention to (pseudo)differential forms on $M$.

Under the Fourier
transform~\eqref{eq:fourier},\eqref{eq:invfourier}, the divergence
operator $\d$ becomes the exterior differential $d$, up to a
multiple of $i$. The multiplication by a multivector field
$H=H(x,x^*)$ becomes the `convolution' (or `cap product'):
\begin{multline*}
    (H*\o)(x,dx)=\const\iint
    Dx^* D(\widebar{dx})\,e^{\i(dx^a-\widebar{dx}^a)x^*_a}\,
    H(x,x^*)\,\o(x,\widebar{dx})\\=
    \int D(\widebar{dx})\,\check
    H(x,dx-\widebar{dx})\,\o(x,\widebar{dx})\,.
\end{multline*}
where $\check H=\check H(x,dx)$ is the inverse Fourier transform of
$H$. In other words, if we denote
\begin{equation*}
    i_H\o:=H*\o\,,
\end{equation*}
we have
\begin{equation}\label{eq:inter}
    i_H=H\Bigl(x,-\i\der{}{dx}\Bigr)\,
\end{equation}
the differential operator, w.r.t. the variables $dx^a$, with the
symbol $H$. It is clear that up to $\i$'s, it is just the classical
internal product of a form by a multivector field, if we deal with
ordinary differential forms and multivectors on an ordinary
manifold.

We immediately get
\begin{thm} On (pseudo)differential forms on $M$, the
Lie derivative w.r.t. the infinitesimal canonical transformation of
$\Pi T^*M$ generated by  $H$ is given by the `Cartan like formula':
\begin{equation}\label{eq:lie}
    \i L_H=[d,i_H],
\end{equation}
where at the r.h.s. stands the commutator of the de Rham
differential  and the interior product by the multivector field $H$
as defined by~\eqref{eq:inter}.
\end{thm}

This is very remarkable. Suppose $M$ is an ordinary manifold. The
operation $i_H$, up to the imaginary unit,  is the familiar interior
product with a multivector field, generalizing the interior product
with a vector field. For Lie derivatives along vector fields one
proves the Cartan formula $L_X=[d,i_X]$. For multivector fields, as
opposed to vector fields, this equation  is taken as the definition
of a `Lie derivative of a differential form along a multivector
field'. In the classical picture it is not seen how these
derivatives corresponds to actual transformations. Now we see that
they are generators of odd canonical transformations acting on
differential forms.

Notice that in general $L_H$ is not a derivation of the algebra
$\O(M)$. Of course,
\begin{equation*}
    L_{\schl H,G \schr}=[L_H,L_G]
\end{equation*}
where at the l.h.s. stands the Schouten bracket of multivector
fields.  Equation~\eqref{eq:lie} implies that the de Rham
differential on $M$ is invariant under the canonical transformations
of $\Pi T^*M$. Again, one can see the $\D$ operator as the de Rham
differential considered together with these extra symmetries.

Some of the arguments of this section were implicit in our earlier
works~\cite{hov:deltabest, hov:max, hov:semi, hov:proclms}.


\begin{thebibliography}{10}

\bibitem{bv:perv}
I.~A. Batalin and G.~A. Vilkovisky.
\newblock Gauge algebra and quantization.
\newblock {\em Phys. Lett.}, 102B:27--31, 1981.

\bibitem{bv:vtor}
I.~A. Batalin and G.~A. Vilkovisky.
\newblock Quantization of gauge theories with linearly dependent generators.
\newblock {\em Phys. Rev.}, D28:2567--2582, 1983.

\bibitem{bv:closure}
I.~A. Batalin and G.~A. Vilkovisky.
\newblock Closure of the gauge algebra, generalized {Lie} equations and
  {Feynman} rules.
\newblock {\em Nucl. Phys.}, B234:106--124, 1984.

\bibitem{hov:deltabest}
O.~M. Khudaverdian {}\footnote{O.~M.~Khudaverdian = H.~M.
Khudaverdian.}.
\newblock Geometry of superspace with even and odd brackets.
\newblock Preprint of the {Geneva University}, {UGVA-DPT} 1989/05-613, 1989.
  Published in: \textit{J. Math. Phys.} 32 (1991), 1934--1937.

\bibitem{hov:bv}
O.~M. Khudaverdian.
\newblock {Batalin-Vilkovisky} formalism and odd symplectic geometry.
\newblock In P.~N. Pyatov and S.~N. Solodukhin, eds., {\em Proceedings of
  the Workshop {``Geometry and Integrable Models''}, {Dubna}, {Russia}, 4-8
  October 1994}. World Scientific Publ., 1995,
\newblock {\tt hep-th 9508174}.


\bibitem{hov:max}
O.~M. Khudaverdian.
\newblock {$\Delta$}-operator on semidensities and integral invariants in the
  {Batalin}--{Vilkovisky} geometry.
\newblock Preprint 1999/135, Max-Planck-Institut f\"ur Mathematik Bonn, 19 p,
  1999,
\newblock \texttt{arXiv:math.DG/9909117}.

\bibitem{hov:semi}
H.~M. Khudaverdian.
\newblock Semidensities on odd symplectic supermanifolds.
\newblock {\em Comm. Math. Phys.}, 247(2):353--390, 2004,
  {\texttt{arXiv:math.DG/0012256}}.

\bibitem{hov:proclms}
H.~M.~ Khudaverdian.
\newblock Laplacians in odd symplectic geometry.
\newblock In Th.~Voronov, ed., {\em Quantization, Poisson Brackets and
  Beyond}, volume 315 of {\em Contemp. Math.}, pages 199--212. Amer. Math.
  Soc., Providence, RI, 2002, \texttt{arXiv:math.DG/0212354}.


\bibitem{tv:laplace1}
H.~M.~Khudaverdian and Th.~Th.~ Voronov.
\newblock On odd {Laplace} operators.
\newblock {\em Lett. Math. Phys.}, 62:127--142, 2002, {\tt
  arXiv:math.DG/0205202}.

\bibitem{tv:laplace2bis}
H.~M.~Khudaverdian and Th.~Th.~ Voronov.
\newblock Geometry of differential operators, and odd {Laplace} operators.
\newblock {\em Russian Math. Surveys}, 58:197--198, 2003,
  \texttt{arXiv:math.DG/0301236}.

\bibitem{tv:laplace2}
H.~M.~Khudaverdian and Th.~Th.~ Voronov.
\newblock On odd {Laplace} operators. {II}.
\newblock In V.~M. Buchstaber and I.~M. Krichever, eds., {\em Geometry,
  Topology and Mathematical Physics. {S.~P.~Novikov's} seminar: 2002--2003},
  volume 212 of {\em Amer. Math. Soc. Transl. (2)}, pages 179--205. Amer. Math.
  Soc., Providence, RI, 2004, \texttt{arXiv:math.DG/0212311}.

\bibitem{yvette:divergence}
Y.~Kosmann-Schwarzbach and J.~Monterde.
\newblock Divergence operators and odd {Poisson} brackets.
\newblock {\em Ann. Inst. Fourier}, 52:419--456, 2002,
  \texttt{arXiv:math.QA/0002209}.


\bibitem{maclane:homology}
S. Mac~Lane.
\newblock {\em Homology}.
\newblock Die Grundlehren der mathematischen Wissenschaften, Bd. 114. Academic
  Press Inc., Publishers, New York, 1963.

\bibitem{ass:bv}
A.~S.~Schwarz.
\newblock Geometry of {B}atalin-{V}ilkovisky quantization.
\newblock {\em Comm. Math. Phys.}, 155(2):249--260, 1993.

\bibitem{ass:symmetry}
A.~S.~Schwarz.
\newblock Symmetry transformations in {B}atalin-{V}ilkovisky formalism.
\newblock {\em Lett. Math. Phys.}, 31(4):299--301, 1994.


\bibitem{severa:originbv}
P. {{\v S}evera}.
\newblock {On the origin of the {BV} operator on odd symplectic supermanifolds},
\texttt{arXiv:math.DG/0506331}.

\end{thebibliography}

\def\cprime{$'$} \def\cprime{$'$}

\end{document}